\documentclass[a4paper]{article}

\usepackage[utf8]{inputenc}
\usepackage[UKenglish]{babel}
\usepackage[T1]{fontenc}           	%Font dei simboli
\usepackage{amsmath}
\usepackage{amssymb}
\usepackage{amsfonts}          		%mathbb{R}
\usepackage[all]{xy}
\usepackage{amscd}
\usepackage{amsthm}
\usepackage{graphics}
\usepackage{color}
\usepackage{mathtools}
\usepackage[a4paper]{geometry}
\usepackage{bbm}
\usepackage{mathrsfs}
\usepackage{lipsum}
\usepackage{indentfirst} 
\usepackage{latexsym} 
\usepackage{xcolor}
\usepackage{cite}
\usepackage[hidelinks]{hyperref}
\usepackage{algorithm}
\usepackage{algpseudocode}
\usepackage{enumitem}
\usepackage{multirow}
\usepackage{multicol}

\newtheorem{thm}{Theorem}[section]

\newtheorem{lem}[thm]{Lemma}

\newtheorem{thm*}{Theorem}%[section*]

\newtheorem{lem*}{Lemma}%[section*]

\newcommand{\R}{\mathbb{R}}

\newcommand{\C}{\mathbb{C}}

\newcommand{\rea}{\textnormal{Re}}

\newcommand{\Pis}{\Pi_{\mathcal{S}}}

\DeclareMathOperator*{\argmin}{arg\,min}

\numberwithin{figure}{section}
\numberwithin{table}{section}

\title{Stabilization of a matrix via a low rank-adaptive ODE}
\author{\scshape{Nicola Guglielmi}\thanks{Gran Sasso Science Institute, L'Aquila, Italy (\href{mailto:nicola.guglielmi@gssi.it}{\tt nicola.guglielmi@gssi.it})}\scshape{ and Stefano Sicilia}\thanks{Gran Sasso Science Institute, L'Aquila, Italy (\href{mailto:stefano.sicilia@gssi.it}{\tt stefano.sicilia@gssi.it})}}
\date{}

\begin{document}
 \maketitle
 
 \vspace{-0.05cm}
 
 \begin{abstract}
  Let $A$ be a square matrix with a given structure (e.g. real matrix, sparsity pattern, Toeplitz structure, etc.) and assume that it is unstable, i.e. at least one of its eigenvalues lies in the complex right half-plane. The problem of stabilizing $A$ consists in the computation of a matrix $B$, whose eigenvalues have negative real part and such that the perturbation $\Delta=B-A$ has minimal norm. The structured stabilization further requires that the perturbation preserves the structural pattern of $A$. We solve this non-convex problem by a two-level procedure which involves the computation of the stationary points of a matrix ODE. We exploit the low rank underlying features of the problem by using an adaptive-rank integrator that follows slavishly the rank of the solution. We show the benefits derived from the low rank setting in several numerical examples, which also allow to deal with high dimensional problems.

 \end{abstract}
 

 \section{Introduction}
 
 Given a matrix $A\in \C^{n\times n}$ with spectrum $\sigma(A)=\left\{ \lambda_1,\dots,\lambda_n\right\}$ ordered such that
 \[
  \rea(\lambda_n) \leq \dots \leq \rea (\lambda_1),
 \]
 we consider the problem of its stabilization, that is we look for a matrix $B$ close to $A$ such that $B$ is a stable matrix, i.e. all its eigenvalues lie in the complex half-plane with real negative part. 
 This problem is well-known in the literature and it has been addressed with different methods (see e.g. \cite{burke2003nonsmooth, gillis2017computing, orbandexivry2013nearest} ). In this work we follow the approach of \cite{guglielmi2017matrix} and we improve it by means of an implementation that highlights the low rank features of the problem. We combine the approach presented in \cite{guglielmi2023rank} with the adaptive-rank integrator of \cite{ceruti2022rank} and we derive a new efficient method. We also introduce an alternative functional that exploits the cubic Hermite interpolating polynomial and we compare it with the one presented in \cite{guglielmi2017matrix}. The results of our new method are similar to the other approaches, but the novelty introduced also allow to deal with higher dimensional problems.
 
 In order to obtain a strict stability, we fix a parameter $\delta>0$ and we further require that the stabilized matrix $B$ is $\delta$-stable, that is all its eigenvalues lie in the set 
 \[
  \C^-_\delta=\{\lambda \in \C : \rea(\lambda)<-\delta\}.
 \]
 Formally, we wish to find a perturbation $\Delta$ of minimum Frobenius norm such that $A+\Delta$ is $\delta$-stable, that is we consider the optimization problem
 \begin{equation}
  \label{prob_opt}
  \argmin_{\sigma(A+\Delta)\subseteq \C^-_\delta} \|\Delta\|_F
 \end{equation}
 and we look for its minimum and its minimizer(s). We denote by $m_\delta(A)$ the number of $\delta$-unstable eigenvalues of $A$ (i.e. the eigenvalues with real part larger than $-\delta$), which is zero if and only if $A$ is $\delta$-stable. For solving problem \eqref{prob_opt}, we use a two-level approach similar to the one presented in \cite{guglielmi2017matrix}. Given a fixed perturbation size $\varepsilon>0$ and a predetermined parameter $\delta>0$ (that ensures strict stability), we rewrite the matrix perturbation $\Delta=\varepsilon E$ with $\|E\|_F=1$ and we minimize in the \textit{inner iteration} the objective functional 
 \begin{equation}
  \label{functional}
  F_\varepsilon(E)=\frac{1}{2}\sum_{i=1}^{n} \left( \left( \rea \left( \lambda_i(A+\varepsilon E)\right)+\delta \right)_+\right)^2=\frac{1}{2}\sum_{i=1}^{m_\delta(A+\varepsilon E)} \left( \rea \left( \lambda_i(A+\varepsilon E)\right)+\delta \right)^2,
 \end{equation}
 where $a_+=\max(a,0)$ is the positive part of $a$. Then the \textit{outer iteration} tunes the perturbation size $\varepsilon>0$ and finds the minimum value $\varepsilon_\star$ such that $F_{\varepsilon_\star}(E(\varepsilon_\star))=0$, where 
 \begin{equation}
  \label{prob_inner}
  E_\star(\varepsilon)=\argmin_{\|E\|_F=1}  F_\varepsilon(E)
 \end{equation}
 is the solution of the optimization problem that arises in the \textit{inner iteration}.
 
 While in \cite{guglielmi2017matrix} the number of addends of the summation in the objective functional is fixed and it is based on an initial guess of the amount of unstable eigenvalues, in our new approach the number of summands depends on $E$ and relies on its current number of $\delta$-unstable eigenvalues.
 This feature makes it possible to exploit the properties of the perturbation $E$ with an adaptive choice of its rank. The optimizers of problem \eqref{prob_opt} are seen as stationary points of a gradient system, which is integrated through a rank-adaptive strategy based on the one presented in \cite{ceruti2022rank}.
 
 The procedure described can be adapted also to the structured stabilization problem. Given $A\in \mathcal{S}$, where $\mathcal{S}\subseteq \C^{n\times n}$ is a subspace of the complex matrices, we consider 
 \begin{equation}
  \label{prob_opt_struct}
  \argmin_{\sigma(A+\Delta)\subseteq \C^-_\delta, \Delta\in \mathcal{S}} \|\Delta\|_F,
 \end{equation}
 that is the structured version of \eqref{prob_opt}. Also in this case we proceed with the two-level approach and we show how the method can be reused in order to exploit all the low rank insights that arise in the unconstrained problem.
 
 The paper is organized as follows. In section 2 we introduce the gradient system for the \textit{inner iteration} of the unstructured problem. In section 3 we show the rank-adaptive integrator used for the solution of the gradient system. In section 4 we investigate the behaviour of the rank during the integration of the gradient system and we show it in several numerical examples. The \textit{outer iteration} is presented in section 5, where an alternative functional for the \textit{inner iteration} is introduced. In section 6 we adapt the procedure to the structured case and in section 7 we show some numerical examples of large dimension.  
 
 \section{Gradient system}
 
 In this section we fix the perturbation size $\varepsilon>0$ and we describe an ordinary differential equation that is used to solve problem \eqref{prob_inner}. In the following, given a matrix whose eigenvalues are simple, we denote by $x_i$ and $y_i$, respectively, the unit left and right eigenvectors associated to $\lambda_i$, such that $x_i^*y_i$ is real and non-negative. We will often exploit the following standard perturbation result for eigenvalues (see \cite{kato2013perturbation}).
 \begin{lem}
  \label{lem_kato}
  Let $\lambda(t)$ be a simple eigenvalue of a differentiable matrix path $A(t)$ in a neighborhood of $t_0$ and let $x(t)$ and $y(t)$ be, respectively, the left and right unit eigenvectors associated. Then $x(t_0)^*y(t_0)\neq 0$ and
   \[
    \dot{\lambda}(t_0)=\frac{x(t_0)^*\dot{A}(t_0)y(t_0)}{x(t_0)^*y(t_0)}.
   \]
 \end{lem}
 We will always suppose that the hypothesis of Lemma \ref{lem_kato} hold true in our setting. This is a generic assumption, since the matrices with at least a $2$-by-$2$ Jordan block form a set of zero measure in $\C^{n\times n}$. 

 In order to find an optimal value of $E$ that minimizes the objective functional $F_\varepsilon(E)$, we introduce a matrix differentiable path $E(t)$ of unit Frobenius norm matrices that depends on a real non-negative time variable $t\geq 0$. By denoting $\mathbb{S}_1$ the unit norm sphere in $\C^{n\times n}$
 \[
  \mathbb{S}_1=\left\{A\in \C^{n\times n}: \|A\|_F=1 \right\},
 \]
 it holds that $E(t)\subseteq \mathbb{S}_1$. In this way it is possible to consider the continuous version $F_\varepsilon(E(t))$ of the objective functional, whose derivative is characterized by the following result.
 
 \begin{lem}
  \label{lem_derf_t}
  Let $E(t)\subseteq \mathbb{S}_1$ be a differentiable path of matrices for $t\in [0,+\infty)$. Let $\varepsilon$ and $\delta$ be fixed. Then $F_\varepsilon(E)$ is differentiable in $[0,+\infty)$ with
  \[
   \frac{d}{dt}F_\varepsilon(E(t))=\varepsilon \rea \langle G_\varepsilon(E(t)), \dot{E}(t)\rangle,
  \]
  where $\langle X,Y\rangle=\textnormal{trace}(X^*Y)$ is the Frobenius inner product and $G_\varepsilon(E(t))$ is the gradient
  \[
   G_\varepsilon(E(t))=\sum_{i=1}^n \gamma_i(t) x_i(t) y_i(t)^*, \qquad
   \gamma_i(t)=\frac{\left(\rea \left(\lambda_i(A+\varepsilon E(t))\right)+\delta\right)_+}{x_i(t)^*y_i(t)}\in \R.
  \]
  \begin{proof}
   We will omit time dependence for brevity. 
   Since 
   \[
    u^*Bv=\langle uv^*,B\rangle, \qquad \forall B\in \C^{n\times n}, \quad \forall u,v\in \C^n
   \]
   and $((x_+)^2)'=2x_+$, Lemma \ref{lem_kato} implies
   \[
    \frac{d}{dt}F_\varepsilon(E)=\sum_{i=0}^n \rea \left( \frac{\varepsilon x_i^*\dot{E} y_i}{x_i^*y_i} \right) \left(\rea\left(\lambda_i(A+\varepsilon E)+\delta \right)\right)_+  
    =\varepsilon \rea \left\langle \sum_{i=1}^n \gamma_i x_i y_i^*, \dot{E}\right\rangle=\varepsilon \rea \langle G_\varepsilon(E), \dot{E}\rangle.
   \]
  \end{proof}
 \end{lem}
 
 The gradient $G=G_\varepsilon(E(t))$ introduced in Lemma \ref{lem_derf_t} is a matrix whose rank depends on the $\delta$-unstable eigenvalues of $A+\varepsilon E(t)$ and it gives the steepest descent direction for minimizing the objective functional, without considering the constraint on the norm of $E$. The next result shows the best direction to follow in order to fulfill the unit norm condition, which is equivalent to $\rea \langle E,\dot{E}\rangle=0$.
 
 \begin{lem}
  \label{lem_opt_grad}
  Given $E\in \mathbb{S}_1$ and $G\in \C^{n\times n}\setminus \{0\}$, the solution of the optimization problem
  \begin{equation}
   \label{optE}
   \argmin_{Z\in \mathbb{S}_1, \  \rea\langle Z,E \rangle=0} \rea \langle G,Z \rangle
  \end{equation}
  is 
  \[
   \alpha Z_\star=-G+\rea \langle G,E\rangle E,
  \]
  where $\alpha>0$ is the normalization parameter.
 \begin{proof}
  Let us consider the vectorized form of the matrices in $\C^{n^2}$. Then the Frobenius product in $\C^{n\times n}$ turns into the standard scalar product of $\C^{n^2}$, which can be seen as a real standard vector product in a $\R^{4n^2}$ and the thesis is straightforward.
 \end{proof}

 \end{lem}
 
 Lemmas \ref{lem_derf_t} and \ref{lem_opt_grad} suggest to consider the matrix ordinary differential equation
 \begin{equation}
  \label{odeE}
  \dot{E}(t)=-G_\varepsilon(E(t))+\rea \langle G_\varepsilon(E(t)),E(t) \rangle E(t),
 \end{equation}
 whose stationary points are zeros of the derivative of the objective functional $F_\varepsilon(E(t))$. Equation \eqref{odeE} is a gradient system for $F_\varepsilon(E(t))$, since along its trajectories
 \[
   \frac{d}{dt}F_\varepsilon(E(t))=\varepsilon\left(-\|G_\varepsilon(E(t))\|_F^2+(\langle G_\varepsilon(E(t)),E(t)\rangle)^2\right)\leq 0
 \]
 by means of the Cauchy-Schwartz inequality, which also implies that the derivative vanishes in $E_\star$ if and only if $E_\star$ is a real multiple of  $G_\varepsilon(E_\star)$. Thanks to the monotonicity property along the trajectories, an integration of this gradient system must lead to a stationary point $E_\star$. The next result provides the important property of the minimizers that reveals the low rank underlying structure of the problem.
 
 \begin{thm}
  \label{thm_prop_stat}
  Let $E_\star$ be a stationary point of \eqref{odeE}, such that $F_\varepsilon(E_\star)>0$. Then $G_\varepsilon(E_\star)\neq 0$ and $E_\star$ is a real multiple of $G_\varepsilon(E_\star)$, that is there exists $\nu\neq 0$ and an integer $1\leq m \leq n$ such that 
  \[
   E_\star=\nu G(E_\star)=\nu \sum_{i=1}^m \gamma_i x_i y_i^*.
  \]
  \begin{proof}
   For all $E\in \mathbb{S}_1$, it holds
   \[
    \rea\langle G_\varepsilon(E),A+\varepsilon E \rangle =\rea \left\langle \sum_{i=1}^m \gamma_i x_i y_i^*,A+\varepsilon E\right\rangle=\rea \left(\sum_{i=1}^m \gamma_i x_i^*(A+\varepsilon E)y_i\right)=
   \]
   \[
    =\sum_{i=1}^m\rea (\gamma_i \lambda_i x_i^*y_i)=\sum_{i=1}^m\rea (\lambda_i) \gamma_i (x_i^*y_i)\geq 0,
   \]
   which means that the gradient $G_\varepsilon(E)$ vanishes if and only if $A+\varepsilon E$ is $\delta$-stable, which is not the case since $F_\varepsilon(E_\star)>0$. Thus, since $E_\star$ is a stationary point, then the right hand side of \eqref{odeE} vanishes and the thesis is straightforward.
  \end{proof}
 \end{thm}

 Theorem \ref{thm_prop_stat} together with the monotonicity property, show that the integration of equation \eqref{odeE} always lead to a low rank stationary point.
 
 \section{A rank-adaptive integrator for the gradient system}
 
 In this section we discuss how to integrate system \eqref{odeE}. Let us assume that for all $t$ the rank $r$ matrix $E(t)\in \C^{n\times n}$ can be decomposed as
 \[
  E(t)=U(t)S(t)V(t)^*, \quad U(t),V(t)\in \C^{n\times r}, \quad S(t)\in \C^{r \times r} \textnormal{ invertible}.
 \] 
 This decomposition generalizes the SVD since it is not required that the matrix $S$ is diagonal. We can rewrite equation \eqref{odeE} as
 \[
  \dot{U}SV^*+U\dot{S}V^*+US \dot{V}^*=-G+\mu US V^*, 
 \]
 where $G=G_\varepsilon(E)$ and $\mu=\rea \langle G,E \rangle$.
 Imposing the gauge conditions $U^*\dot{U}=V^*\dot{V}=0$ yields the system
 \begin{equation}
  \label{odeUSV}
  \begin{dcases}
   \dot{U}=-(I-UU^*)GV S^{-1}\\
   \dot{S}=-U^*GV+\mu S \\
   \dot{V}=(I-VV^*)G^*US^{-*}
  \end{dcases},
 \end{equation}
 which is equivalent to \eqref{odeE}. The structure of the matrix $S$ is not diagonal in general, but the decomposition assumed for $E$ holds along all the trajectory. System \eqref{odeUSV} consists of two matrix ODEs of dimension $n$-by-$r$ and one of dimension $r$-by-$r$, where $r$ is the rank of $E(t)$. 
 
 A classical integration, e.g. by means of Euler method, of system \eqref{odeUSV} is not suitable. Indeed the presence of $S^{-1}$ may cause numerical issues and moreover this integration does not capture the change of the rank of $G$ along the trajectory. To overcome this problem, we exploit a rank-adaptive strategy similar to the one exposed in \cite{ceruti2022rank}. We define
 \[
  f(E):=-G_\varepsilon(E)+\rea\langle G_\varepsilon(E),E\rangle E
 \]
 as the left hand side of \eqref{odeE} such that $\dot{E}=f(E)$ and we fix a tolerance $\tau$.
 To update the perturbation path $E(t)=U(t)S(t)V(t)^*$ from $t_0$ to $t_1$, we start from $E_0=E(t_0)=U(t_0)S(t_0)V(t_0)^*=U_0S_0V_0^*$ of rank $r_0$ and we get $E_1=E(t_1)=U(t_1)S(t_1)V(t_1)^*=U_1S_1V_1^*$ of rank $r_1$ by performing the following steps: 
 \begin{enumerate}
  \item Set $\rho=\min(2r_0,n)$.
  \item Compute augmented basis $\widehat{U}\in \C^{n\times \rho}$ and $\widehat{V}\in \C^{n\times \rho}$
  
  \smallskip
  
  \textbf{K-step}: Integrate from $t_0$ to $t_1$ the $n$-by-$r_0$ differential equation
  \[
   K(t_0)=U_0S_0, \qquad \dot{K}(t)=f\left(K(t)V_0^*\right)V_0.
  \]
  Perform a QR factorization of $(K(t_1),U_0)$, save the first $\rho$ columns in $\widehat{U}\in\C^{n\times \rho}$ and compute $\widehat{M}=\widehat{U}^*U_0\in \C^{\rho\times r_0}$.
  
  \smallskip
  
  \textbf{L-step}: Integrate from $t_0$ to $t_1$ the $n$-by-$r_0$ differential equation
  \[
   L(t_0)=V_0S_0^*, \qquad \dot{L}(t)=f\left(U_0L(t)^*\right)^*U_0.
  \]
  Perform a QR factorization of $(L(t_1),V_0)$, save the first $\rho$ columns in $\widehat{V}\in\C^{n\times \rho}$ and compute $\widehat{N}=\widehat{V}^*V_0\in \C^{\rho\times r_0}$.

  \item Augment and update $S$
  
  \smallskip
  
  \textbf{S-step}: Integrate from $t_0$ to $t_1$ the $\rho$-by-$\rho$ differential equation
  \[
   \widehat{S}(t_0)=\widehat{M}S_0\widehat{N}^*, \qquad \dot{\widehat{S}}(t)=\widehat{U}^*f\left(\widehat{U}\widehat{S}(t)\widehat{V}^*\right)\widehat{V}.
  \] 
  
  \item Adapt the rank for the updated matrices
  
  \smallskip
  
  \textbf{Truncation}: Compute the SVD $\widehat{S}(t_1)=\widehat{P}\widehat{\Sigma}\widehat{Q}^*$ and choose the new rank $r_1\leq \rho$ such that
  \[
   \left(\sum_{i=r_1+1}^\rho \sigma^2_i\right)^{\frac{1}{2}}\leq \tau, \qquad \widehat{\Sigma}=\textnormal{diag}(\sigma_i).
  \]
  Define $S_1$ as the $r_1$-by-$r_1$ diagonal main sub-matrix of $\widehat{\Sigma}$ and denote by $P_1,Q_1\in \C^{\rho\times r_1}$ the first $r_1$ columns of $\widehat{P}$ and $\widehat{Q}$ respectively.
  
  \item Return $U_1=\widehat{U}P_1 \in \C^{n\times r_1}$, $V_1=\widehat{V}Q_1 \in \C^{n\times r_1}$ and $S_1\in \C^{r_1\times r_1}$.

 \end{enumerate}
 In \cite{ceruti2022rank} it is shown that this algorithm computes an approximation of $U(t),S(t)$ and $V(t)$ with an error proportional to the tolerance $\tau$ and the time step $t_1-t_0$. Moreover this algorithm allows the truncation of the rank, according to the tolerance $\tau$ in order to adapt the size of the invertible matrix $S$ along the trajectory.
 
 \section{Rank adaptivity for fixed perturbation size}
 
 In this section we show with three different examples how the integrator introduced in the previous section chooses the adaptive rank of the perturbation. Unless otherwise stated, we set the parameter $\delta=10^{-3}$.
 
 \subsection{An illustrative example}
 
 Consider the matrix
 \begin{equation}
  \label{mat_ill}
  A=\left(\begin{array}{rrrrrrrrrr}
     0 & 1 & 1 & 1 & -1 & 0 & -1 & 0 & 0 & 0 \\
     1 & -1 & 0 & 1 & 1 & 0 & 1 & 0 & 0 & 0 \\
     -1 & 0 & -1 & -1  & -1 & 1 & 1 & 1 & 0 & 0  \\
     1 & 0 & 0 & -1 & 1 & -1 & -1 & 1 & 0 & 0 \\
     0 & 0 & -1 & 1 & 0 & 1 & 1 & -1 & 0 & 0 \\
     0 & -1 & 1 & 1 & -1 & 0 & 0 & 1 & 1 & 0 \\
     -1 & 1 & -1 & 1 & 1 & 0 & -1 & 0 & 1 & 1 \\
     0 & 0 & 1 & -1 & -1 & 1 & 1 & 1 & -1 & 1 \\
     0 & 0 & 0 & 0 & 0 & 0 & 0 & -1 & 1 & -1 \\
     0 & 0 & 0 & 0 & 0 & 0 & 0 & 0 & -1 & 1 \\
    \end{array}\right)\in \C^{10 \times 10},
 \end{equation}
 which has 6 unstable eigenvalues. Table \ref{tab_illu} contains the results of the functional $F_\varepsilon(E_\star)$ and the maximum rank of $E(t)$ achieved during the \textit{inner iteration} with different choices of $\varepsilon$ and of the tolerance $\tau$ In particular $\varepsilon=2.38$ is close to be $\varepsilon_\star$.
 
 \begin{table}
  \begin{center}
   \begin{tabular}{||c||cc||cc||}
    \hline
    \multirow{2}{*}{$\tau$} & \multicolumn{2}{c||}{$\varepsilon=2$}                     & \multicolumn{2}{c||}{$\varepsilon=2.38$}                  \\ \cline{2-5} 
                        & \multicolumn{1}{c|}{$F_\varepsilon(E_\star)$} & max rank & \multicolumn{1}{c|}{$F_\varepsilon(E_\star)$} & max rank \\ \hline
    $10^{-1}$               & \multicolumn{1}{c|}{4.1060}                   & 6        & \multicolumn{1}{c|}{2.1548}                   & 6        \\ \hline
    $10^{-2}$               & \multicolumn{1}{c|}{1.9701}                   & 6        & \multicolumn{1}{c|}{0.2932}                   & 6        \\ \hline
    $10^{-3}$               & \multicolumn{1}{c|}{1.7830}                   & 7        & \multicolumn{1}{c|}{0.0687}                   & 7        \\ \hline
    $10^{-4}$               & \multicolumn{1}{c|}{0.7325}                   & 8        & \multicolumn{1}{c|}{0.0058}         & 9        \\ \hline
    $10^{-5}$               & \multicolumn{1}{c|}{0.5289}                   & 9        & \multicolumn{1}{c|}{$8\cdot 10^{-7}$}         & 9        \\ \hline
    $10^{-6}$               & \multicolumn{1}{c|}{0.5483}                   & 9        & \multicolumn{1}{c|}{$2\cdot 10^{-7}$}         & 9        \\ \hline
   \end{tabular}
   \caption{Illustrative matrix \eqref{mat_ill}: numerical results of the \textit{\textit{inner iteration}}}
   \label{tab_illu}
  \end{center}
\end{table}
 
 In both cases it is possible to observe how much the value of the objective function decreases when the tolerance is lowered, which is explained by the higher adaptive rank of the perturbation. This means that lower values of the tolerance lead to more accurate results, but they also increase the rank of the perturbation and hence the computational cost. Hence a suitable choice of $\tau$ is crucial to balance this two factors.

 \subsection{Grcar matrix}
 
 The rank adaptive procedure is not effective in all cases. For instance let us consider the Grcar $n$-by-$n$ matrix, with $n\geq 5$, that is an Hessenberg and Toeplitz matrix of the form
 \begin{equation}
  \label{mat_grcar}
  G_n=\left(\begin{array}{cccccccccc}
     1 & 1 & 1 & 1 & 0 & & \\
     -1 & 1 & 1 & 1 & 1 & 0  \\
     0 & \ddots & \ddots & \ddots & \ddots & \ddots & \ddots  \\
     & \ddots & \ddots & \ddots & \ddots & \ddots & \ddots & 0 \\
     & & 0 & -1 & 1 & 1 & 1 & 1 \\
    \end{array}\right)\in \C^{n \times n}.
 \end{equation}
 The eigenvalues of $G_n$ are all unstable and also quite sensitive. We consider $n=20$ and we show the results of the  \textit{inner iteration} with different choices of $\varepsilon$ and of the tolerance $\tau$.
 
 \begin{table}
  \begin{center}
   \begin{tabular}{||c||cc||cc||}
    \hline
    \multirow{2}{*}{$\tau$} & \multicolumn{2}{c||}{$\varepsilon=4$}                     & \multicolumn{2}{c||}{$\varepsilon=5.5$}                   \\ \cline{2-5} 
                        & \multicolumn{1}{c|}{$F_\varepsilon(E_\star)$} & max rank & \multicolumn{1}{c|}{$F_\varepsilon(E_\star)$} & max rank \\ \hline
    $10^{-1}$               & \multicolumn{1}{c|}{24.7887}                  & 20       & \multicolumn{1}{c|}{6.1836}                   & 20       \\ \hline
    $10^{-2}$               & \multicolumn{1}{c|}{15.9710}                  & 20       & \multicolumn{1}{c|}{0.7401}                   & 20       \\ \hline
    $10^{-3}$               & \multicolumn{1}{c|}{13.2096}                  & 20       & \multicolumn{1}{c|}{0.2218}                   & 20       \\ \hline
    $10^{-4}$               & \multicolumn{1}{c|}{11.9657}                  & 20       & \multicolumn{1}{c|}{0.0082}                   & 20       \\ \hline
    $10^{-5}$               & \multicolumn{1}{c|}{10.8504}                  & 20       & \multicolumn{1}{c|}{0.0001}                   & 20       \\ \hline
    $10^{-6}$               & \multicolumn{1}{c|}{10.6692}                  & 20       & \multicolumn{1}{c|}{$7\cdot 10^{-7}$}         & 20       \\ \hline
   \end{tabular}
   \caption{Grcar matrix \eqref{mat_grcar}: numerical results of the \textit{\textit{inner iteration}}}
   \label{tab_grcar}
  \end{center}
 \end{table}
 
 In this case the fact that all the eigenvalues are unstable provides a perturbation of full (or almost full) rank and the adaptive rank integrator seems not to be effective. During the integration, some eigenvalues may become stable before the other ones and after that they begin to cross back and forth the imaginary axis. Consequently they activate or disable their respective rank-1 component in the gradient, producing a minimal change of the rank (e.g. 18 or 19) that is not worth to exploit with the adaptive integrator.  

 \subsection{Smoke matrix} 
 
 In this example we show more in detail the changes of the rank perturbation during the integration. Let $S_n$ be the $n\times n$ Smoke matrix from the \texttt{gallery} function of \texttt{Matlab}
 \begin{equation}
  \label{mat_smoke}
  S_n=\left(\begin{matrix}
    \zeta_1 & 1 & \\
    & \zeta_2 & 1 & \\
    & & \ddots & \ddots \\
    & & & \ddots & 1\\
    & & & & \zeta_{n-1} & 1 \\
    1 & & & & & \zeta_n
  \end{matrix}\right),
 \end{equation}
 whose diagonal contains the $n$ distinct roots of $1$ 
 \[
  \zeta_j=\mathrm{e}^{\frac{2\pi \mathrm{i}j}{n}}, \qquad j=1,\dots,n
 \]
 such that $\zeta_j^n=1$ for all $j$. The characteristic polynomial $p(\lambda)=\det(S_n-\lambda I_n)$ associated to $S_n$ is
 \[
  p(\lambda)=(-1)^{n+1}+\prod_{j=1}^n (\zeta_j-\lambda)=(-1)^n\left(-1+\prod_{j=1}^n (\lambda-\zeta_j)\right)=(-1)^n(\lambda^n-2)
 \]
 and hence the eigenvalues of $S$ are equally distributed along the circle of radius $\sqrt[n]{2}$. For $n$ even the matrix $S_n$ has half eigenvalues stable and half unstable. We consider $n=20$ and we collect the results of the \textit{inner iteration} in Table \ref{tab_smoke1}.
 
 \begin{table}
  \begin{center}
   \begin{tabular}{||c||cc||}
    \hline
    \multirow{2}{*}{$\tau$} & \multicolumn{2}{c||}{$\varepsilon=2.5$}                   \\ \cline{2-3} 
                        & \multicolumn{1}{c|}{$F_\varepsilon(E_\star)$} & max rank \\ \hline
    $10^{-1}$               & \multicolumn{1}{c|}{28.1862}                  & 11       \\ \hline
    $10^{-2}$               & \multicolumn{1}{c|}{0.2931}                   & 11       \\ \hline
    $10^{-3}$               & \multicolumn{1}{c|}{0.0136}                   & 13       \\ \hline
    $10^{-4}$               & \multicolumn{1}{c|}{0.0013}                   & 15       \\ \hline
    $10^{-5}$               & \multicolumn{1}{c|}{$8\cdot 10^{-5}$}                   & 16       \\ \hline
    $10^{-6}$               & \multicolumn{1}{c|}{$3\cdot 10^{-6}$}         & 17       \\ \hline
   \end{tabular}
   \caption{Smoke matrix \eqref{mat_smoke}: numerical results of the \textit{\textit{inner iteration}} }
   \label{tab_smoke1}
  \end{center}
\end{table}
 
\begin{figure}
  \begin{center}
   \includegraphics[scale=0.54]{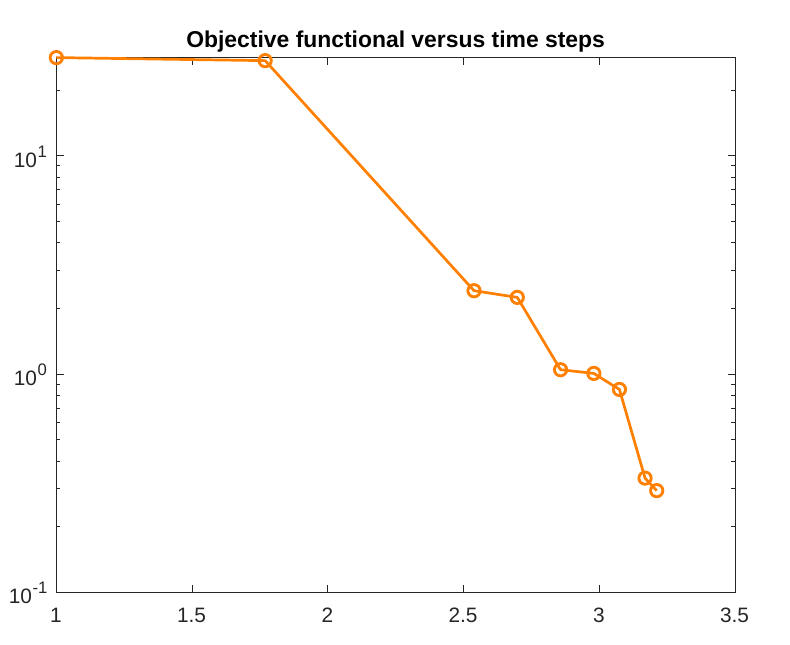}
   \includegraphics[scale=0.54]{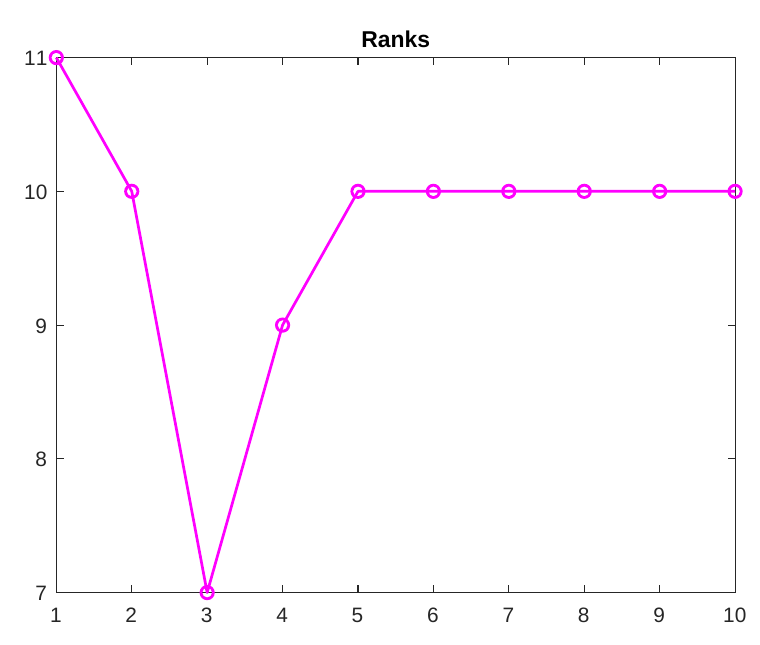}
   \includegraphics[scale=0.54]{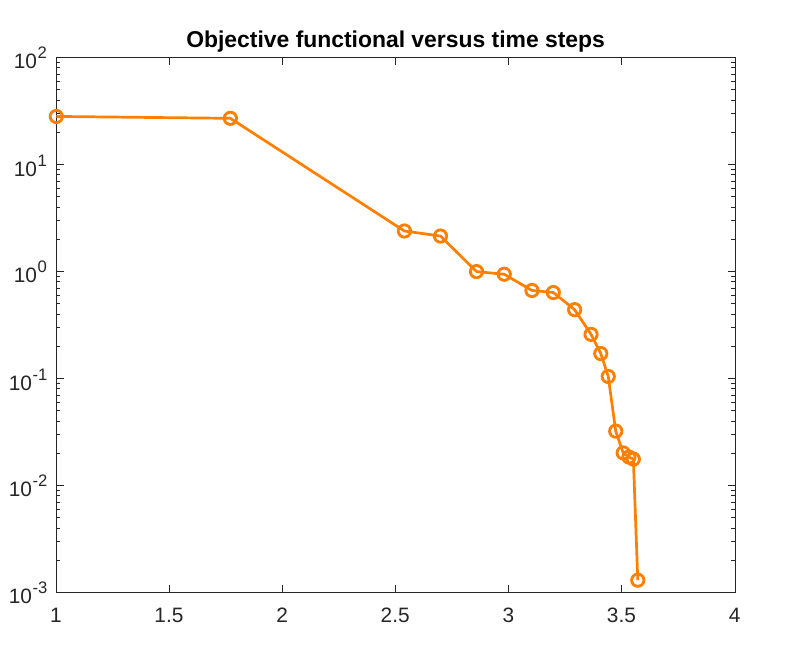}
   \includegraphics[scale=0.54]{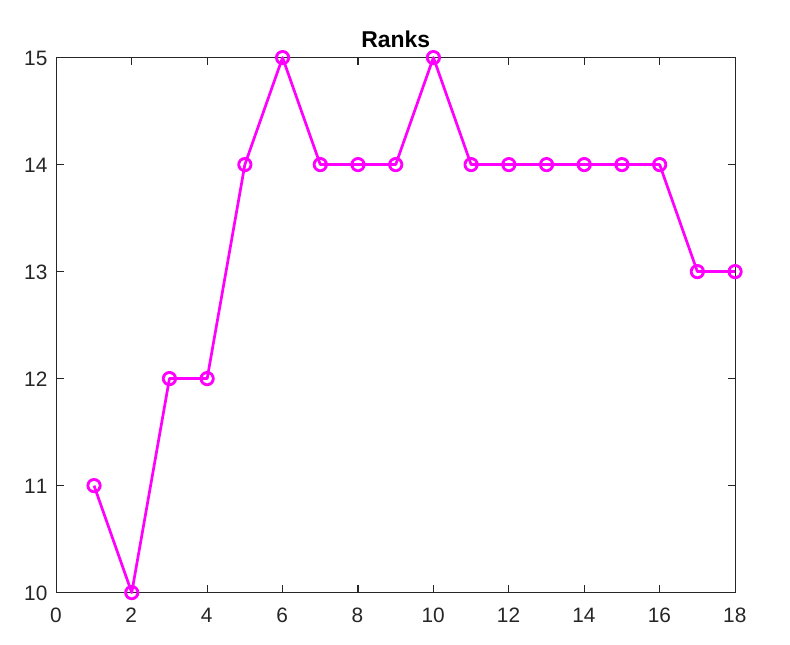}
   \includegraphics[scale=0.54]{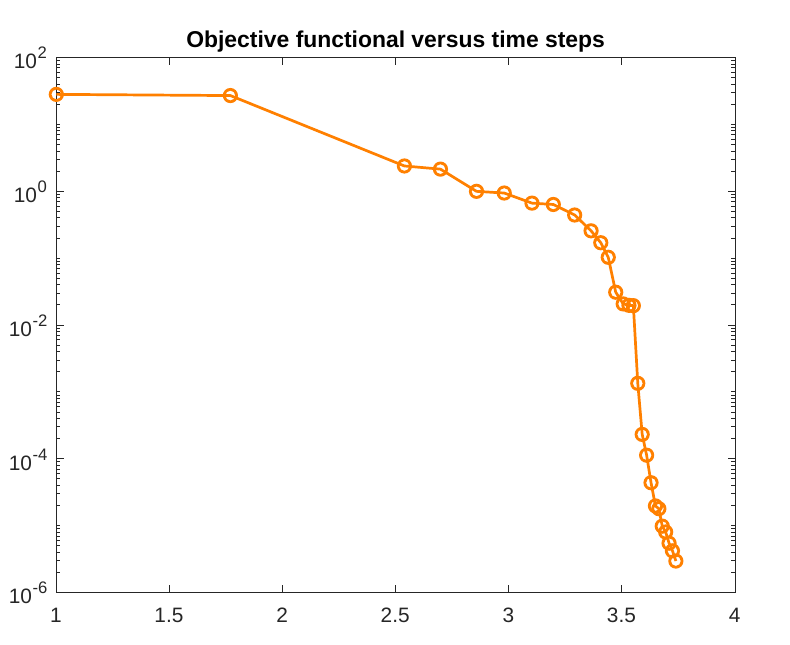}
   \includegraphics[scale=0.54]{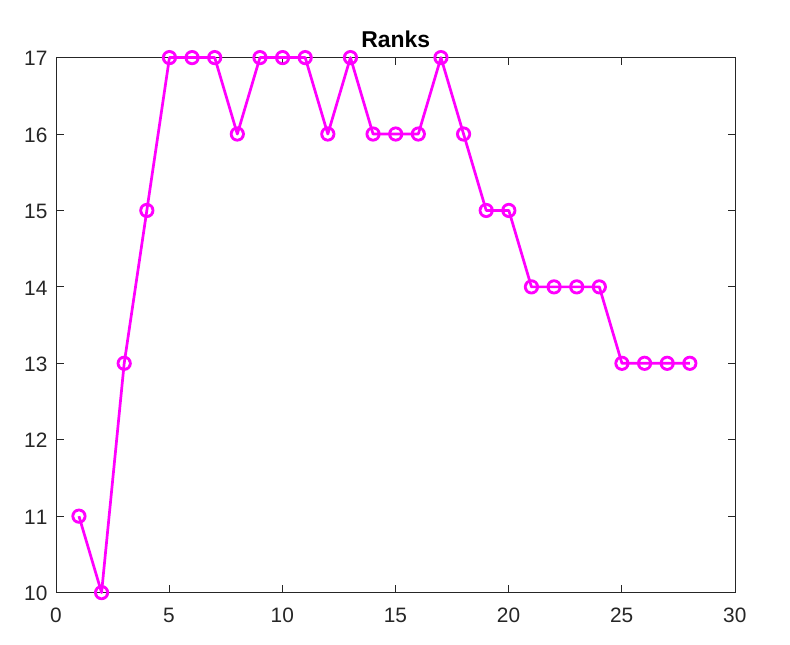}
  \end{center}
  \caption{Smoke matrix: functional and ranks in the \textit{\textit{inner iteration}} for $\tau=10^{-2}$ (up), $\tau=10^{-4}$ (middle) and $\tau=10^{-6}$ (down)}
  \label{fig_smoke}
 \end{figure}

 Figure \ref{fig_smoke} shows the trend of the rank and the objective functional. We can notice how the rank is increased or decreased by the integrator, but definitely it stabilizes when it approaches the stationary point.

 \section{Outer iteration: tuning the perturbation size}
  
 Once that a computation of the optimizers is available for a given $\varepsilon>0$, we need to determine an optimal value for the perturbation size $\varepsilon_\star$.
 Let $E_\star(\varepsilon)$ be a solution of the optimization problem \eqref{prob_inner} and consider the function
 \[
  \varphi(\varepsilon)=F_\varepsilon(E_\star(\varepsilon)).
 \]
 This function is non-negative and we define $\varepsilon_\star$ as the smallest zero of $\varphi$. Assuming that the unstable eigenvalues of $A+\varepsilon E_\star(\varepsilon)$ are simple, for $0\leq \varepsilon < \varepsilon_\star$, yields that $\varphi$ is a smooth function in the interval $[0,\varepsilon_\star)$.
 The aim of the \textit{\textit{outer iteration}} is to approximate $\varepsilon_\star$, which is the solution of the optimization problem \eqref{prob_opt}. In order to solve this problem, we use a combination of the well-known Newton and bisection methods, which provides an approach similar to \cite{guglielmi2015low,guglielmi2016method} or \cite{guglielmi2023rank}. If the current approximation $\varepsilon$ lies in a left neighborhood of $\varepsilon_\star$, it is possible to exploit Newton's method, since $\varphi$ is smooth there; otherwise in a right neighborhood we use the bisection method. The following result provides a simple formula for the first derivative of $\varphi$ which is cheap to compute making the Newton's method easy to apply.
 
 \begin{lem}
  \label{lem_derfeps}
  For $0\leq \varepsilon<\varepsilon_\star$ it holds
  \[
   \varphi'(\varepsilon)=\frac{d}{d\varepsilon} F_\varepsilon(E_\star(\varepsilon))=\langle G_\varepsilon(E_\star(\varepsilon)),E_\star(\varepsilon)\rangle=-\|G_\varepsilon(E_\star)\|_F\leq 0.
  \]
  \begin{proof}
   As shown for the time derivative formula, we get
   \[
    \varphi'(\varepsilon)=\frac{d}{d\varepsilon} \left(\frac{1}{2}\sum_{i=1}^{n} \left( \left( \rea \left( \lambda_i(A+\varepsilon E(\varepsilon))\right)+\delta \right)_+\right)^2\right)=
   \]
   \[
    =\sum_{i=1}^{n}\left(\rea \left( \lambda_i(A+\varepsilon E(\varepsilon))\right)+\delta \right)_+\frac{d}{d\varepsilon} \left( \rea \left( \lambda_i(A+\varepsilon E(\varepsilon)) \right)_+\right)=\rea \langle G_\varepsilon(E_\star(\varepsilon)),E_\star(\varepsilon)+\varepsilon E'_\star(\varepsilon)\rangle.
   \]
   where $E'_\star(\varepsilon)$ is the derivative with respect to $\varepsilon$ of $E_\star(\varepsilon)$. 
   Since $E_\star$ is a unit norm stationary point of \eqref{odeE} and a zero of the derivative of the objective functional $F_\varepsilon$, then $G_\varepsilon(E_\star(\varepsilon))$ is a negative multiple of $E_\star$. Thus 
   \[
    G_\varepsilon(E_\star)=-\|G_\varepsilon(E_\star)\|_F E_\star(\varepsilon), \qquad
    \rea \langle G_\varepsilon(E_\star(\varepsilon)), E'_\star(\varepsilon)\rangle=\frac{\rea \langle E_\star(\varepsilon), E'_\star(\varepsilon) \rangle}{\rea \langle G_\varepsilon(E_\star(\varepsilon)), E_\star(\varepsilon)\rangle}=0.
   \]

  \end{proof}
 \end{lem}

 \subsection{Another functional}
 
 For a deeper analysis of the method developed, we introduced also the alternative objective functional $\Phi_\varepsilon$. Given $0<\delta_1<\delta_2$ we define it as
 \[
  \Phi_\varepsilon(E)=\frac{1}{2}\sum_{i=1}^n \psi(\rho_i)\left(\rho_i+\delta_2\right)_+^2, \qquad \rho_i=\rea(\lambda_i(A+\varepsilon E)),
 \]
 where $\psi\in C^1(\R)$ is the cubic Hermite interpolating polynomial
 \[
  \psi(x)=\begin{cases}
   0 & x<-\delta_2 \\
   \frac{(x+\delta_2)^2(2x+3\delta_1-\delta_2)}{(\delta_1-\delta_2)^3} & -\delta_2\leq x \leq -\delta_1 \\
   1 & x>-\delta_1
  \end{cases}, 
  \qquad
  \psi'(x)=\begin{cases}
   \frac{6(x+\delta_2)(x+\delta_1)}{(\delta_1-\delta_2)^3} & -\delta_2\leq x \leq -\delta_1 \\
   0 & \textnormal{otherwise}
  \end{cases}
 \]
 such that
 \[
  \psi(\delta_1)=1, \qquad \psi(\delta_2)=\psi'(\delta_1)=\psi'(\delta_2)=0.
 \]
 The aim of this functional is to soften the passage of an eigenvalue from unstable to stable and viceversa, by the introduction of the Hermite interpolating polynomial $\psi(x)$. Lemma \ref{lem_derf2_t} shows that the previous theory can be applied also in this case, by means of a slight change of the gradient associated to the objective functional. In particular the new gradient is made up by the same rank-one perturbations of the previous version, but the coefficients associated are different. 
 
 \begin{lem}
  \label{lem_derf2_t}
  Let $E(t)\subseteq \mathbb{S}_1$ be a differentiable path of matrices for $t\in [0,+\infty)$. Let $\varepsilon$ and $\delta$ be fixed. Then $\Phi_\varepsilon(E)$ is differentiable in $[0,+\infty)$ with
  \[
   \frac{d}{dt}\Phi_\varepsilon(E(t))=\varepsilon \rea \langle \Gamma_\varepsilon(E(t)), \dot{E}(t)\rangle,
  \]
  where $\Gamma_\varepsilon(E(t))$ is the gradient of $\Phi_\varepsilon$
  \[
   \Gamma_\varepsilon(E(t))=\sum_{i=1}^n \kappa_i(t) x_i(t) y_i(t)^*, \qquad
   \kappa_i(t)=\frac{(\rho_i+\delta_2)_+\left(\psi'(\rho_i)(\rho_i+\delta_2)_++2\psi(\rho_i)\right)}{2x_i(t)^*y_i(t)}\geq 0,
  \]
  with $\rho_i(t)=\rea \left(\lambda_i(A+\varepsilon E(t))\right)$.
  \begin{proof}
   The proof is the same as that of Lemma \ref{lem_derf_t}, since for any $i=1,\dots,n$
   \[
    \frac{d}{dt}\left(\psi(\rho_i)(\rho_i+\delta_2)_+^2\right)=\dot{\rho}_i\psi'(\rho_i)(\rho_i+\delta_2)_+^2+2\psi(\rho_i)\dot{\rho}_i(\rho_i+\delta_2)_+=
   \]
   \[
    =\dot{\rho}_i(\rho_i+\delta_2)_+\left(\psi'(\rho_i)(\rho_i+\delta_2)_++2\psi(\rho_i)\right)=\langle \kappa_i x_i y_i^*,\dot{E} \rangle.
   \]

  \end{proof}
 \end{lem}
 Also the $\textit{outer iteration}$ can be adapted to the functional $\Phi_\varepsilon$. In practice we have always considered $\delta_2=2\delta_1$ and $\delta=\delta_1=10^{-3}$, but other choices are possible.
 
 \subsection{Smoke matrix}
 
 Let us consider again the Smoke matrix $S\in \R^{n\times n}$. For $n=30$ we generate a random orthogonal matrix $[Q,\sim]=\texttt{qr}\left(\texttt{randn}(n)+\mathrm{i}\cdot \texttt{randn}(n)\right)$ (we set \texttt{rng(1)}) and we applied the algorithm on the Smoke-like matrix 
 \begin{equation}
  \label{mat_smokelike}
  A=QSQ^*.
 \end{equation}
 In table \ref{tab_smoke_comp} we consider both the functionals $F_\varepsilon$ and $\Phi_\varepsilon$ and we report the results for different integration strategies: we study the difference between the adaptive rank approach and the fixed rank method proposed in  proposed in \cite{guglielmi2017matrix} with several values of the rank $r$. In all the experiments we set
 \[
  \tau_{\textnormal{inner}}=10^{-9}, \qquad \tau_{\textnormal{outer}}=10^{-9}, \qquad \tau_{\textnormal{rank}}=10^{-8}, \qquad  \textnormal{maxit}_{\textnormal{inner}}=250, \qquad \textnormal{maxit}_{\textnormal{outer}}=300,
 \]
 where $\tau$ stands for tolerance, maxit for the maximum number of iteration allowed and inner and outer refers to the \textit{inner} and \textit{outer iterations} respectively.   
 
 \begin{table}
  \begin{center}
  \begin{tabular}{|c|c|c|c|c|c|}
   \hline
                                             & Rank feature    & $\varepsilon_\star$ & rank$(E_\star(\varepsilon_\star))$ & Functional & CPU time (sec) \\ \hline \hline
   \multirow{10}{*}{$F_\varepsilon$}            & Adaptive        & 3.7547              & 20                                 & $6.1879\cdot 10^{-11}$ & 10.8693            \\ \cline{2-6} 
                                             & Fixed (15)      & 3.8083              & 15                                 & $1.4060\cdot 10^{-11}$ & 12.0071            \\ \cline{2-6} 
                                             & Fixed (16)      & 3.9802              & 16                                 & $1.8237\cdot 10^{-11}$ & 12.5968            \\ \cline{2-6} 
                                             & Fixed (17)      & 3.7500              & 17                                 & $2.2933\cdot 10^{-11}$ & 10.3501            \\ \cline{2-6} 
                                             & Fixed (18)      & 3.7783              & 18                                 & $4.0462\cdot 10^{-11}$ & 9.0029             \\ \cline{2-6} 
                                             & Fixed (19)      & 3.6478              & 19                                 & $1.1513\cdot 10^{-11}$ & 11.5401            \\ \cline{2-6} 
                                             & Fixed (20)      & 3.9539              & 20                                 & $8.1013\cdot 10^{-12}$ & 11.9893            \\ \cline{2-6} 
                                             & Fixed (21)      & 3.8667              & 21                                 & $4.5873\cdot 10^{-11}$ & 16.5965            \\ \cline{2-6} 
                                             & Fixed (22)      & 3.8818              & 22                                 & $1.1716\cdot 10^{-11}$ & 10.9894            \\ \cline{2-6} 
                                             & Fixed (23)      & 3.9601              & 23                                 & $1.9844\cdot 10^{-11}$ & 10.0265            \\ \hline \hline
   \multirow{10}{*}{$\Phi_\varepsilon$}         & Adaptive        & 3.7104              & 20                                 & $1.0000\cdot 10^{-9}$ & 16.9300            \\ \cline{2-6} 
                                             & Fixed (15)      & 3.6475              & 15                                 & $1.0000\cdot 10^{-9}$ & 20.7475            \\ \cline{2-6} 
                                             & Fixed (16)      & 3.7491              & 16                                 & $1.0000\cdot 10^{-9}$ & 12.4132            \\ \cline{2-6} 
                                             & Fixed (17)      & 3.7723              & 17                                 & $1.0000\cdot 10^{-9}$ & 13.0365            \\ \cline{2-6} 
                                             & Fixed (18)      & 3.6416              & 18                                 & $2.2043\cdot 10^{-9}$ & 25.3195            \\ \cline{2-6} 
                                             & Fixed (19)      & 3.6949              & 19                                 & $5.4971\cdot 10^{-10}$ & 13.3895            \\ \cline{2-6} 
                                             & Fixed (20)      & 3.8774              & 20                                 & $5.6181\cdot 10^{-9}$ & 25.1754            \\ \cline{2-6} 
                                             & Fixed (21)      & 4.0715              & 21                                 & $9.8411\cdot 10^{-10}$ & 10.0267            \\ \cline{2-6} 
                                             & Fixed (22)      & 3.9193              & 22                                 & $1.0000\cdot 10^{-9}$ & 24.9207            \\ \cline{2-6} 
                                             & Fixed (23) & 3.8251              & 23                                 & $1.0000\cdot 10^{-9}$ & 16.3523            \\ \hline \hline
   \multirow{3}{*}{G-S algorithm} & BCD             & 3.7087              & 24                                 & $5.0000\cdot 10^{-7}$ & 0.6573                 \\ \cline{2-6} 
                                             & Grad            & 3.3706              & 28                                 & $5.0000\cdot 10^{-7}$ & 0.3317                 \\ \cline{2-6} 
                                             & FGM             & 3.0975              & 30                                 & $1.0000\cdot 10^{-6}$ &  0.2548                  \\ \hline
   \end{tabular}
   \label{tab_smoke_comp}
   \caption{Comparison with the adaptive rank and fixed rank integrators of the Smoke-like matrix \eqref{mat_smokelike}. The results of the algorithm of Gillis and Sharma (G-S) are also reported with three different methods.}   
   \end{center}
 \end{table}
 
 For the standard functional $F_\varepsilon$, the adaptive-rank integrator provides the third smallest distance $\varepsilon_\star=3.7547$, close to the best one provided by the fixed rank method with $r=19$ that is $\varepsilon_\star=3.6478$. A similar analysis also holds for the Hermite functional case $\Phi_\varepsilon$. Even though the adaptive rank does not provide the best result, it detects quite well the rank of the best optimizer found by the fixed rank procedure. In this way it is possible to avoid to perform several times the fixed rank integration in order to improve the result, which is more expensive especially when it is not available a GPU.
  
 \begin{figure}
  \begin{center}
   \includegraphics[scale=0.52]{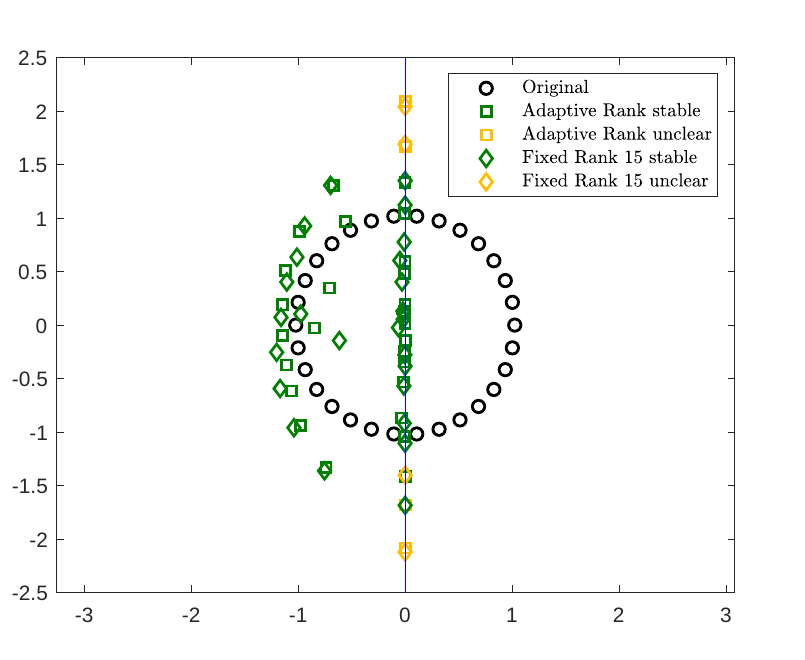}
   \includegraphics[scale=0.52]{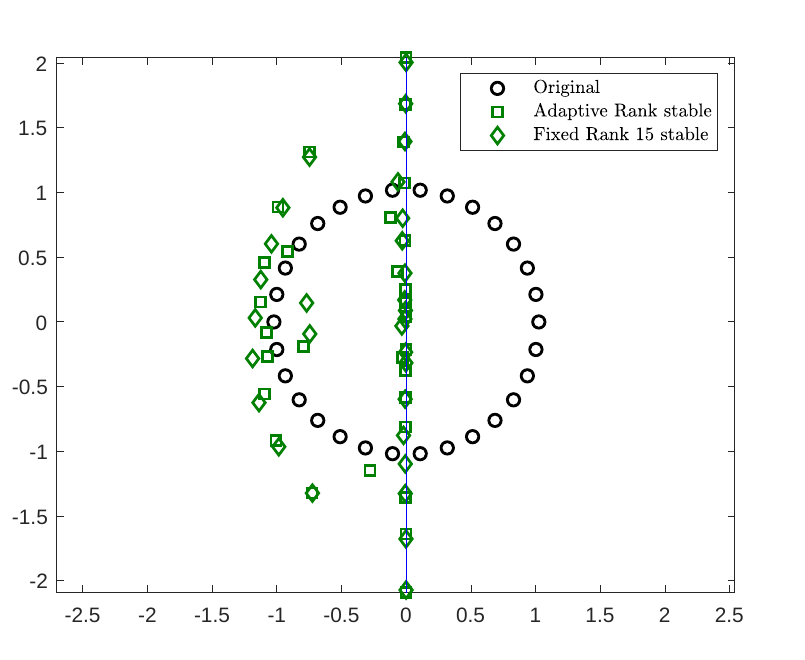}
  \end{center}
  \caption{Smoke-like matrix \eqref{mat_smokelike}: original eigenvalues (black circles) and stabilized ones (an eigenvalue $\lambda$ is green if  $\rea(\lambda)<-\delta$, while orange if $-\delta\leq \rea(\lambda)\leq 0$). On the left the functional considered is  $F_\varepsilon$, while on the right $\Phi_\varepsilon$. }
  \label{fig_eig_smoke}
 \end{figure}
 
 In all the experiments we observed an alignment along the imaginary axis of the unstable eigenvalues, as shown in Figure \ref{fig_eig_smoke}. This means that all the algorithms move the unstable eigenvalues towards the imaginary axis, but they also change some of the stable ones making them closer to have zero real part. This is the reason behind the fact that the optimal fixed rank (19 in this example) is not the original number of unstable eigenvalues (15).
 
 We also compared the results of our method with the algorithms proposed by Gillis and Sharma in \cite{gillis2017computing} and two of them (Grad and FGM) provide smaller distances in less computational time. However this type of algorithm does not take into account the rank properties of the problem and hence the perturbations found have higher rank, almost full. This means that for low dimensional example the algorithm of Gillis and Sharma is preferable, but in higher dimension it could be too expensive to perform. Moreover the method of \cite{gillis2017computing} is not easy to extend to the structured version of the problem.

 \section{Structured distance via a low rank adaptive ODE}
 
 Let now $\mathcal{S}\in\C^{n\times n}$ be a linear subspace, such as a prescribed sparsity pattern, Toeplitz matrices, Hamiltonian matrices, etc.. Given an unstable matrix $A\in \mathcal{S}$, we look for a  matrix $\Delta\in \mathcal{S}$ with smallest norm that stabilizes $A$. Formally we want to solve the optimization problem
 \[
  \argmin_{\Delta\in \mathcal{S}} \left\{\|\Delta\|_F : \sigma(A+\Delta)\subseteq \C^-_{\delta}\right\},
 \]
 which is a generalization of \eqref{prob_opt}. For the solution of this structured optimization problem we follow the same approach used in the previous sections. Let $\Pis$ be the orthogonal projection, with respect to the inner Frobenius product, to the subspace $\mathcal{S}$. Generally it is easy to compute an explicit formula for $\Pis$ and some examples can be found in \cite{guglielmi2023rank}. By proceeding as in the unstructured case, we project onto $\mathcal{S}$ the gradient $G_\varepsilon$ (see for instance \cite{guglielmi2017matrix}) to get the new system
 \begin{equation}
  \label{SodeE}
  \dot{E}=-\Pis(G_\varepsilon(E))+\real\langle \Pis(G_\varepsilon(E)),E\rangle E.
 \end{equation}
 Equation \eqref{SodeE} represents a gradient system where the gradient $G_\varepsilon$ has been replaced by $\Pis G_\varepsilon$. All the results for the unconstrained case extend to the structured system, with the replacement of the structured gradient. However $\Pis(G_\varepsilon(E))$ is generally not low rank, since the property of $G_\varepsilon$ is now subject to the presence of the projection $\Pis$ and cannot be exploited as before. But it turns out that also in this case we can formulate a low rank adaptive ODE that takes into account also the structure constraint.
 Let us assume that $E=\Pis Y$, for a certain $Y$ that has the same rank $r$ of $G_\varepsilon(E)$. As in \cite{guglielmi2023rank}, let us consider the ODE
 \begin{equation}
  \label{odeY}
  \dot{Y}=-P_Y G_\varepsilon(\Pis Y) +\real\langle P_YG_\varepsilon(\Pis Y), \Pis Y\rangle Y,
 \end{equation}
 where $P_Y$ is the projection onto the tangent space $\mathcal{T}_Y\mathcal{M}_r$ at the rank $r$ manifold $\mathcal{M}_r$. An explicit expression for $P_Y$ is given by (see \cite{guglielmi2017matrix})
 \[
  P_Y(A)=A-(I-UU^*)A(I-VV^*)
 \]
 where $Y=USV^*$ is an SVD decomposition of $Y$ (here $S$ is required to be invertible, but it may not be diagonal). If the rank of $Y$ is fixed, this is a low rank ODE whose stationary points are explicitly related to the ones of the gradient system and in particular this correspondence is bijective. Equation \eqref{odeY} is not a gradient system, but it can be proved, similarly as in \cite{guglielmi2023rank} and \cite{guglielmi2023low}, that its integration locally converges to a stationary point $Y_\star$ that represents a stationary point $E_\star=\Pis Y_\star$ of the original gradient system \eqref{SodeE}. Indeed when the trajectory is close to a stationary point, the rank of the gradient stabilizes and thus it is possible to proceed as in \cite[Theorem 4.4]{guglielmi2023rank}. Thus it is possible to solve \eqref{odeY} by means of the adaptive-rank integrator, so that we can exploit the low rank insights of this problem even in the structured case.

 \section{Numerical examples for the structured case}
 
 In this section we investigate the behaviour of the algorithm for the structured optimization problem in some numerical examples, including matrices with large dimension. In all the cases we consider as structure the sparsity pattern of the matrix. The code used is available upon request.
  
 \subsection{Pentadiagonal Toeplitz matrix}
 
\begin{figure}[H]
  \begin{center}
   \includegraphics[scale=0.6]{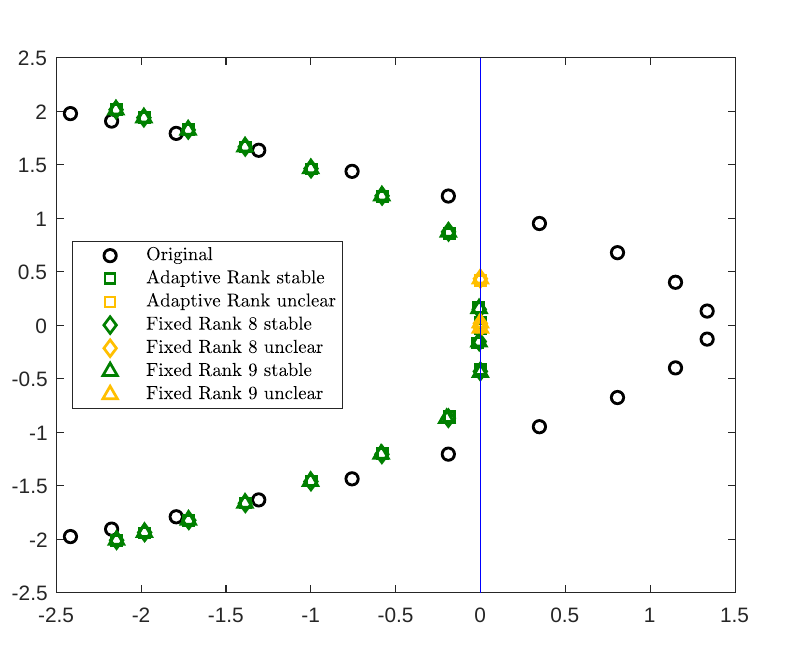}
  \end{center}
  \caption{Pentadiagonal Toeplitz matrix: original eigenvalues (black circles) and stabilized ones (an eigenvalue $\lambda$ is green if  $\rea(\lambda)<-\delta$, while orange if $-\delta\leq \rea(\lambda)\leq 0$).}
  \label{fig_eig_pentadiag}
 \end{figure}

 Let us consider the pentadiagonal Toeplitz matrix
 \begin{equation}
  \label{mat_pentadiag}
  P=\left(\begin{matrix}
    -\frac{1}{2} & 1 & 1 &  \\
    1 & -\frac{1}{2} & 1 & 1 & \\
    1 & 1 & \ddots & \ddots & \ddots \\
    & \ddots & \ddots & \ddots & \ddots & 1\\
    & & \ddots & \ddots & -\frac{1}{2} & 1 \\
    & & & 1 & 1 & -\frac{1}{2} \\
  \end{matrix}\right)\in \R^{20\times 20}.
 \end{equation}
 We apply the method where the solution of the \textit{inner iteration} is obtained by the integration of \eqref{odeY} and we collect the results in  Table \ref{tab_pentadiagS} and Figure \ref{fig_eig_pentadiag}.
 
 \begin{table}
  \begin{center}
  \begin{tabular}{|c|c|c|c|c|}
  \hline
  & $\varepsilon_\star$ & rank$(E_\star(\varepsilon_\star))$ & $F(E_\star(\varepsilon_\star))$ & CPU time (seconds) \\ \hline \hline
  Adaptive rank                         & 2.9024              & 10                                 & $1.0000\cdot 10^{-9}$             & 6.9174             \\ \hline
  Fixed rank (8)                        & 2.9033              & 8                                  & $4.7522\cdot 10^{-10}$             & 9.8106             \\ \hline
  Fixed rank (9)                        & 2.9093              & 9                                  & $1.9073\cdot 10^{-10}$             & 10.6489            \\ \hline
  Fixed rank (10)                       & 2.9011              & 10                                 & $8.0005\cdot 10^{-10}$             & 9.4543             \\ \hline
  \multicolumn{1}{|l|}{Fixed rank (11)} & 2.9023              & 11                                 & $8.4959\cdot 10^{-10}$             & 9.3636             \\ \hline
  \end{tabular}
  \end{center}
  \caption{Pentadiagonal Toeplitz matrix \eqref{mat_pentadiag} results with functional $F_\varepsilon$.}
  \label{tab_pentadiagS}
 \end{table}
  
 All the algorithms provide a similar distance, but in this case it seems that the rank-adaptive integrator is quicker than the fixed rank methods. The results are identical also when the algorithms are applied for the functional $\Phi$.
 
 \subsection{Brusselator matrix}
 
 Now we consider the Brusselator matrix\footnote{See https://math.nist.gov/MatrixMarket/data/NEP/brussel/rdb800l.html} from the NEP collection. This matrix has size $n=800$ with non-zeros $nnz=4640\approx 5.8 n$ and it arises from a two-dimensional reaction-diffusion model in chemical engineering. It has two conjugate unstable eigenvalues that are close to the imaginary axis and some eigenvalues with negative real part close to zero. We applied the algorithm with $\delta=10^{-3}$ and parameters
 \[
  \tau_{\textnormal{inner}}=10^{-9}, \qquad \tau_{\textnormal{outer}}=10^{-9}, \qquad \tau_{\textnormal{rank}}=10^{-8}, \qquad  \textnormal{maxit}_{\textnormal{inner}}=250, \qquad \textnormal{maxit}_{\textnormal{inner}}=300.
 \]
 For the fixed rank method, we were only able to select $r=2$, since in the other cases the computations in the \texttt{Matlab} function \texttt{eigs} did not converge.
 
 \begin{table}
  \begin{center}
   \begin{tabular}{|c|c|c|c|c|c|}
    \hline 
                                            & Rank feature & $\varepsilon_\star$ & rank$(E_\star(\varepsilon_\star))$ & Functional             & CPU time (seconds) \\ \hline \hline
    \multirow{4}{*}{$F_{\varepsilon_\star}$}    & Adaptive     & 0.9458              & 4                                  & $5.1538\cdot 10^{-10}$ & 171           \\ \cline{2-6} 
                                            & Fixed (2)    & 0.9431              & 2                                  & $5.0269\cdot 10^{-10}$ & 182           \\ \cline{2-6} 
                                            & Fixed (3)    & -                   & -                                  & -                      & -                  \\ \cline{2-6} 
                                            & Fixed (4)    & -                   & -                                  & -                      & -                  \\ \hline \hline
    \multirow{4}{*}{$\Phi_{\varepsilon_\star}$} & Adaptive     & 0.9437              & 4                                  & $5.2447\cdot 10^{-10}$ & 372           \\ \cline{2-6} 
                                            & Fixed (2)    & 0.9374              & 2                                  & $4.4979\cdot 10^{-11}$ & 808           \\ \cline{2-6} 
                                            & Fixed (3)    & -                   & -                                  & -                      & -                  \\ \cline{2-6} 
                                            & Fixed (4)    & -                   & -                                  & -                      & -                  \\ \hline
   \end{tabular}
   \caption{Brusselator matrix: features of the solution of the \textit{\textit{outer iteration}}}
   \label{tab_brus}
  \end{center}
 \end{table}
 
 The results show that for the functional $F$ the algorithm is quicker even if provides slightly worse distances than those associated to $\Phi$. In this case it seems that the adaptive rank captures the fact that two stable eigenvalues are close to the imaginary axis and considered them in the gradient, so that they do not become unstable. Hence the final adaptive rank is 4 instead of 2.
 
 \begin{figure}
  \begin{center}
   \includegraphics[scale=0.52]{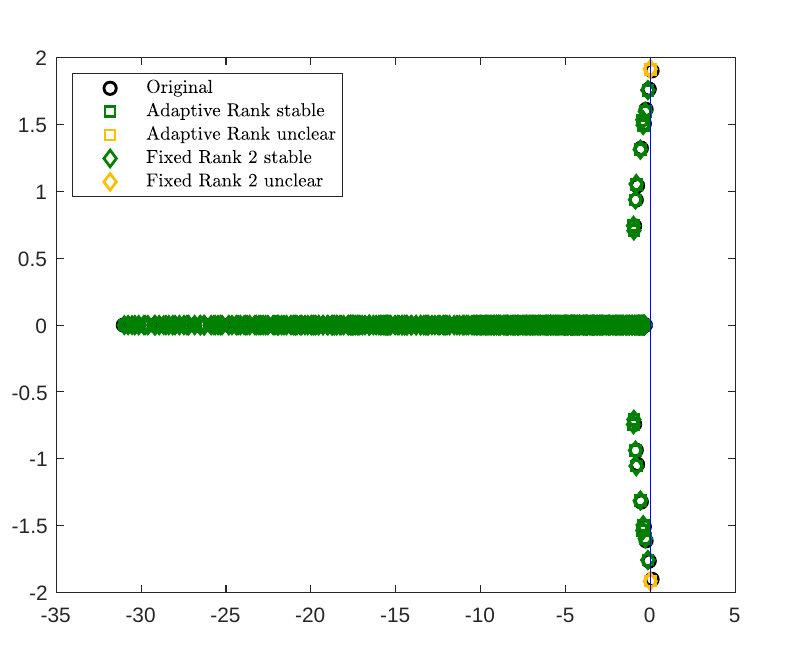}
   \includegraphics[scale=0.52]{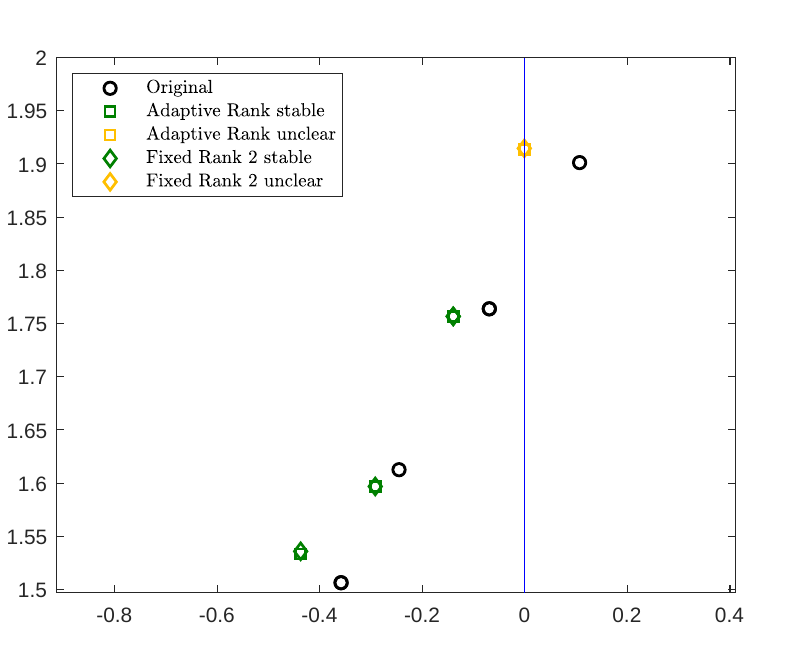}
  \end{center}
  \caption{Brusselator matrix: original eigenvalues (black circles) and stabilized ones (an eigenvalue $\lambda$ is green if  $\rea(\lambda)<-\delta$, while orange if $-\delta\leq \rea(\lambda)\leq 0$). On the right a zoom of the image of all the eigenvalues in the left.}
  \label{fig_eig_brus}
 \end{figure}
 
 \subsection{Fidap matrix}
  
\begin{figure}
  \begin{center}
   \includegraphics[scale=0.52]{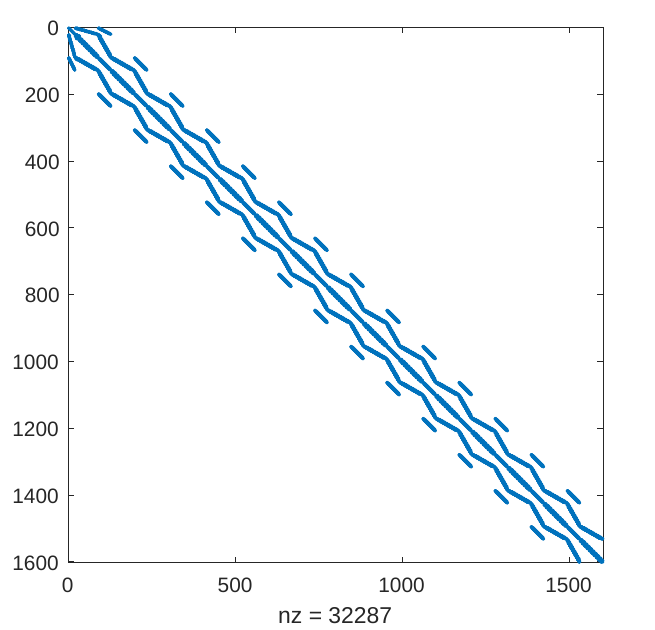}
   \includegraphics[scale=0.52]{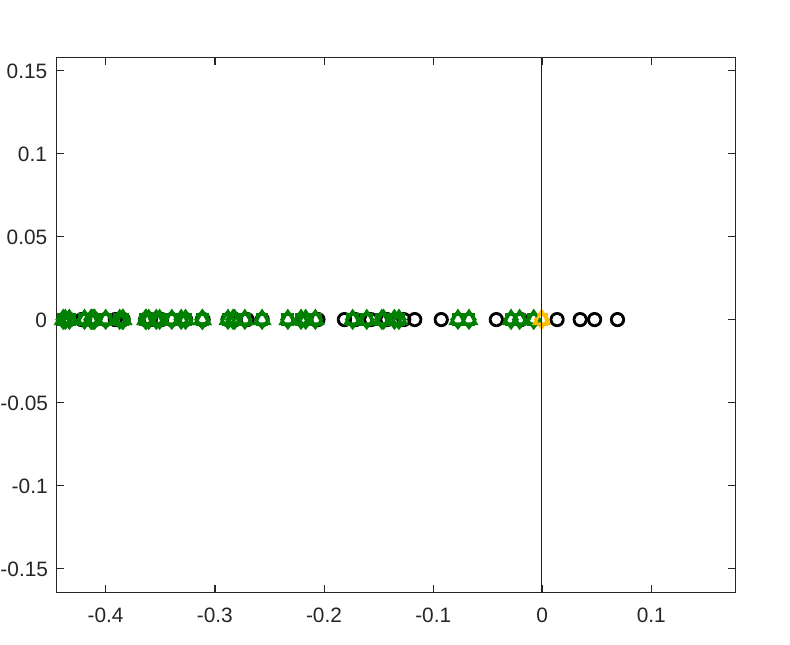}
  \end{center}
  \caption{Fidap matrix: on the left its structural pattern, on the right its original eigenvalues (black circles) and the stabilized ones (an eigenvalue $\lambda$ is green if  $\rea(\lambda)<-\delta$, while orange if $-\delta\leq \rea(\lambda)\leq 0$). On the right a zoom of the image of all the eigenvalues in the left.}
  \label{fig_fidap}
 \end{figure}

 Finally we consider the Fidap matrix from the SPARSKIT collection\footnote{See https://math.nist.gov/MatrixMarket/data/SPARSKIT/fidap/fidap004.html}, which arises in fluid dynamics modeling. This matrix has size $1601\times 1601$ and it is symmetric. In our example we consider the shifted matrix $A-\frac{3}{2}I$ so that the number of unstable eigenvalues reduces to 4. In this case the parameters chosen are
 \[
  \tau_{\textnormal{inner}}=10^{-9}, \qquad \tau_{\textnormal{outer}}=10^{-9}, \qquad \tau_{\textnormal{rank}}=10^{-8}, \qquad  \textnormal{maxit}_{\textnormal{inner}}=150, \qquad \textnormal{maxit}_{\textnormal{inner}}=200.
 \]
 Table \ref{tab_fidap} and Figure \ref{fig_fidap} show the results of the algorithm applied on the shifted matrix.
 
 \begin{table}
  \begin{center}
   \begin{tabular}{|c|c|c|c|c|c|}
    \hline
                                            & Rank feature & $\varepsilon_\star$ & rank$(E_\star(\varepsilon_\star))$ & Functional        & CPU time (seconds) \\ \hline \hline
    \multirow{4}{*}{$F_{\varepsilon_\star}$}    & Adaptive     & 0.1820              & 6                                  & $1.0000\cdot 10^{-9}$  & 60                 \\ \cline{2-6} 
                                            & Fixed (4)    & 0.1821              & 4                                  & $9.9999\cdot 10^{-10}$  & 50                 \\ \cline{2-6} 
                                            & Fixed (5)    & 0.1820              & 5                                  & $1.0000\cdot 10^{-9}$  & 55                 \\ \cline{2-6} 
                                            & Fixed (6)    & 0.1820              & 6                                  & $1.0000\cdot 10^{-9}$  & 49                 \\ \hline \hline
    \multirow{4}{*}{$\Phi_{\varepsilon_\star}$} & Adaptive     & 0.1819              & 4                                  & $9.3534\cdot 10^{-10}$ & 135                \\ \cline{2-6} 
                                            & Fixed (4)    & 0.1819              & 4                                  & $7.3800\cdot 10^{-10}$ & 285                \\ \cline{2-6} 
                                            & Fixed (5)    & 0.1819              & 5                                  & $7.2320\cdot 10^{-10}$ & 371                \\ \cline{2-6} 
                                            & Fixed (6)    & 0.1819              & 6                                  & $7.2490\cdot 10^{-10}$ & 350                \\ \hline
   \end{tabular}
   \caption{Fidap matrix: results with functional $F_\varepsilon$.}
   \label{tab_fidap}
  \end{center}
 \end{table}
 
  Also in this case we observe that the distance provided by the functional $\Phi$ is slightly lower than the one associated to the functional $F$, but this improvement is not significant enough to justify the higher computational time needed by $\Phi$ with respect to the $F$.

 \section*{Acknowledgments}
  Nicola Guglielmi acknowledges that his research was supported by funds from the Italian MUR (Ministero dell'Universit\`{a} e della Ricerca) within the PRIN 2021 Project ``Advanced numerical methods for time dependent parametric partial differential equations with applications'' and the Pro3 Project ``Calcolo scientifico per le scienze naturali, sociali e applicazioni: sviluppo metodologico e tecnologico''. Nicola Guglielmi and Stefano Sicilia are affiliated to the Italian INdAM-GNCS (Gruppo Nazionale di Calcolo Scientifico).

 \nocite{*}
 \bibliographystyle{plain}
 \bibliography{biblio}

\end{document}